\documentclass[10pt]{amsart}

\usepackage{amsfonts,amssymb,amsmath,amsthm}
\usepackage{url}
\usepackage{enumerate}
\usepackage{xcolor}
\usepackage{hyperref}

\newtheorem{thrm}{Theorem}[section]
\newtheorem{lem}[thrm]{Lemma}
\newtheorem{prop}[thrm]{Proposition}
\newtheorem{cor}[thrm]{Corollary}
\newtheorem{definition}[thrm]{Definition}

\numberwithin{equation}{section}

\newcommand{\Z}{\mathbb{Z}}
\newcommand{\F}{\mathbb{F}}

\newcommand{\CC}{\mathcal C}
\newcommand{\vr}{\vspace*{0.5cm}}

\title{ LCD and Self Orthogonal twisted group codes over finite commutative chain rings}

\author{Samir Assuena $^{1.}$ } 
\address{1. Centro Universit\'ario da FEI}
\email{samir.assuena@fei.edu.br}

\author{Andr\'e Luiz Martins Pereira $^{2.}$}
\address{2. Universidade Federal Rural do Rio de Janeiro .} 
\email{almp1980@ufrrj.br}


\begin{document}

\begin{abstract}
In this paper, we shall study  $k$-Galois LCD constacyclic group codes over finite commutative chain rings with identity. In particular, we shall characterize Galois LCD constacyclic codes over finite commutative chain ring with identity in terms of its idempotent generators and the classical involution using the twisted group ring structures and find some good LCD codes. 
\end{abstract} 
\maketitle

\section{Introduction}

Linear codes with complementary duals (abbreviated LCD) are linear codes whose intersection with their dual is trivial. When they are binary, they play an important role in armoring implementations against side-channel attacks and fault injection attacks. 

Linear complementary dual codes have importance in data storage, communications systems and security too.

These codes have been studied for improving the security of information on sensitive devices against side-channel attacks (SCA) and fault non-invasive attacks, see \cite{CC}, and have found use in data storage and communications systems.

Carlet and Guilley, in \cite{CG}, also investigated the application of binary LCD codes against side-channel attacks (SCA) and fault tolerant injection attacks (FIA). Also, in \cite{TD}, the authors constructed explicity LCD codes and have explicit efficient encoding and decoding algorithms .

In \cite{FZ}, Fan and Zhang, introduced the concept of \textit{k-Galois form}, which is a generalization of Euclidean and Hermitian inner products and Liu, Fan and Liu, in \cite{LFL}, studied $k$-Galois LCD codes. Also, in \cite{ZL}, Liu introduced the {\textit{generalized Galois inner product}} for commutative rings.


So, this paper is devoted to classify the $k$-Galois LCD constacyclic codes over finite commutative chain rings in terms of idempotents in twisted group algebras of cyclic groups. 


Let $R$ be a finite commutative ring with identity, $\CC$ be a linear code over $R^{n}$, that is, $\CC$ is a $R$-submodule of $R^{n}$ and let $\lambda$ be an element of ${\mathcal{U}}(R)$, the group of units of $R$. We say that $\CC$ is a \textit{$\lambda$-constacyclic code} if 

\begin{center}
$(c_{0}, c_{1}, \cdots, c_{n-1}) \in \CC \Longrightarrow (\lambda c_{n-1}, c_{0}, \cdots, c_{n-2}) \in \CC $
\end{center}

\noindent for all $(c_{0}, c_{1}, \cdots, c_{n-1}) \in \CC$.

When $\lambda=1$, we have so called \textit{cyclic codes} and, when $\lambda=-1$, we have \textit{negacyclic codes}. Thus, constacyclic codes are generalization of cyclic and negacyclic codes and they have been studied for many authors (\cite{BR}, \cite{CDLW}, \cite{D}). Also, constacyclic codes can be realized as ideals in polynomial factor ring $\displaystyle\frac{R[x]}{\left\langle x^{n}-\lambda\right\rangle}$. 

Given $x=(x_{0}, x_{1}, \cdots, x_{n-1})$ and $y=(y_{0}, y_{1}, \cdots, y_{n-1})$ two elements of a linear code  $\CC$, the \textit{Hamming distance} between $x$ and $y$ is the number

\begin{center}
	$d_{H}(x,y)=|\{i: \, x_{i}\neq y_{i}, \,\, 0\leq i\leq n-1\}|$.
\end{center}

\noindent and the \textit{weight} of $x$ is 

\begin{center}
	$w_{H}(x)=d(x,0)=|\{i: \, x_{i}\neq 0, \,\, 0\leq i\leq n-1\}|$.
\end{center}

It is well-known that, for a linear code $\CC$, we have $d_{H}(x,y)=w_{H}(x-y)$, for all $x,y \in \CC$.


Let $G$ be a group and $A$ be an abelian group. A map

\begin{center}
$\alpha: G\times G \longrightarrow A$
\end{center}

\noindent is a 2-\textit{cocycle} if , for all $x,\,y$ and $z$ in $G$, we have

\begin{center}
$\alpha(x,y)\alpha(xy,z)=\alpha(y,z)\alpha(x,yz)$.
\end{center}

\noindent and a map $t:G\times G \longrightarrow A$ is a 2-\textit{coboundary} if there is a map $\delta:G \longrightarrow A $ such that

\begin{center}
$t(x,y)=\delta(x)\delta(y)\delta(xy)^{-1}$.
\end{center}

As usual, the set of all 2-cocycles will be denoted by $Z^{2}(G,A)$ and the set of all 2-coboundary will be denoted by $B^{2}(G,A)$. Finally, we say that a 2-cocycle $\alpha$ is \textit{normalized} if $\alpha(x,1)=\alpha(1,x)=\alpha(1,1)=1$, for all $x \in G$. Notice that, since $\alpha$ is a 2-cocycle, we can replace $\alpha$ by $\alpha'$ given by 

\begin{center}
$\alpha'(x,y)=\displaystyle\frac{\alpha(x,y)}{\alpha(1,1)}$
\end{center}

\noindent which is a normalized 2-cocycle.  From now on, we assume that all 2-cocycles are normalized.

Let $R$ be a commutative ring and $G$ be a group. The \textit{twisted group ring $R^{\gamma}G$ of G over} $R$ is the associative $R$-algebra with basis $\overline{G}=\{\overline{g}, \, g\in G\}$, which is a copy of $G$, and the multiplication is defined on the basis as

\begin{center}
$\overline{g}\cdot \overline{h}=\gamma(g,h)\overline{gh}$
\end{center}

\noindent where $\gamma(g,h)$ is an element of ${\mathcal{U}}(R)$, the group of units of $R$.

The mapping $\gamma:G\times G \longrightarrow {\mathcal{U}}(R)$ is called \textit{twisting} and there are many different possibilities for $R^{\gamma}G$ depending on the choice of the twisting. For instance, the group ring $RG$ of $G$ over $R$ is a twisted group ring with $\gamma(g,h)=1$. Furthermore, the associative condition on the multiplication implies that

\begin{center}
$\gamma(g,h)\gamma(gh,k)=\gamma(h,k)\gamma(g,hk)$
\end{center}

\noindent and, for this reason, $\gamma$ is a 2-cocycle. 

When $G=C_{n}=\left\langle g\right\rangle$, a cyclic group of order $n$ and $R$ a commutative, we have the following

 \begin{thrm}\label{lem1}  {\rm{\cite [Theorem 2.3 ]{ARP}}} 
 	Let $C_n = \langle g \rangle$ be a cyclic group of order $n$ and R be a commutative ring. Then, the twisted group ring $R^{\gamma}C_{n}$ is isomorphic to $R^{\gamma_{\lambda}}C_{n}$ where
	
	\begin{equation*}
 		\gamma_\lambda(g^i, g^j) = \begin{cases}
 			1, \quad i+j<n\\
 			\lambda, \quad i+j \ge n
 		\end{cases}
 	\end{equation*}

 \end{thrm}

It is possible make a \textit{diagonal} change of basis by replacing each $\overline{g}$ by $\widetilde{g}=\delta(g)\overline{g}$ for some $\delta(g) \in {\mathcal{U}}(R)$ and, with this change of basis, $R^{\gamma}G$ is realized in a second way as a twisted group ring of $G$ over $R$ with twisting

\begin{center}
$\widetilde{\gamma}(g,h)=\delta(g)\delta(h)\delta(gh)^{-1}\gamma(g,h)$.
\end{center}

In this case, we say that $\gamma$ and $\widetilde{\gamma}$ are \textit{cohomologous}.

\begin{lem} \label{1.1} {\rm{\cite [Lemma 2.1]{PA}}}
The following relations hold in $R^{\gamma}G$

\begin{enumerate}
	\item [i.] $1=\gamma(1,1)^{-1}\overline{1}$
	\item [ii.] For all $g \in G$,

\begin{center}
$\overline{g}^{-1}=\gamma(g,g^{-1})^{-1}\gamma(1,1)^{-1}\overline{g^{-1}}=\gamma(g^{-1},g)^{-1}\gamma(1,1)^{-1}\overline{g^{-1}}$
\end{center}
	
\end{enumerate}

\end{lem}

Let $C_{n}=\left\langle g \mid g^{n}=1\right\rangle$ be a cyclic group of order $n$, $R$ be a finite commutative ring and $R^{\gamma}C_{n}$ the twisted group algebra with 

\begin{center}
		$\gamma_{\lambda}(g^{j},g^{k})= \left\{
		\begin{array}{lll}
			\lambda,    & {\rm{if}} \,\, j+k\geq n\\
			1, &  {\rm{if}} \,\, j+k < n
		\end{array}
		\right.$
\end{center}

\noindent where $\lambda$ is a unity of $R$. Thus, $\overline{g}^{2}=\overline{g}\cdot\overline{g}=\gamma(g,g)\overline{g^{2}}$, so we can make a diagonal change of basis and replace $\overline{g^{k}}$ by $\overline{g}^{k}$, for all $k, \,\, 1\leq k\leq n$. Thus, there exists a unity $a \in R$  such that $\overline{g}^{n}=a\cdot 1$ which implies that $R^{\gamma}C_{n}$ is a commutative ring.

In \cite{SA}, Assuena gave an explicit proof that $\lambda$-constacyclic codes can be also realized as ideals in twisted group algebras $R^{\gamma_{\lambda}}C_{n}$ of cyclic group over a finite commutative ring where

\begin{center}
		$\gamma_{\lambda}(g^{j},g^{k})= \left\{
		\begin{array}{lll}
			\lambda,    & {\rm{if}} \,\, j+k\geq n\\
			1, &  {\rm{if}} \,\, j+k < n.
		\end{array}
		\right.$
\end{center}

\noindent for some $\lambda \in {\mathcal{U}}(R)$. Thus, we say that any ideal $\CC$ of a twisted group algebra $R^{\gamma}G$ of a finite group $G$ over a finite ring $R$ is a \textit{twisted group code}.

In this paper, we shall consider twisted group rings $R^{\gamma_{\lambda}}C_{n}$ of commutative local rings $R$ over finite cyclic groups $C_{n}$ such that $char(R/J(R))$ does not divide $n$, where $J(R)$ denotes the Jacobson radical of $R$.

\section{Some basic results}

In this section, we shall give definitions and some known results which have elementary proofs in twisted group algebras language.

Let $R$ be a finite commutative ring with identity with $p^{n}$ elements, where $p$ is a rational prime, $G$ be a finite group and $R^{\gamma}G$ the twisted group ring of $G$ over $R$. Let $\sigma$ be an automorphism of $R$ with order $m$. Given $\alpha=\displaystyle\sum_{g \in G}\alpha_{g}\overline{g}$, $\beta=\displaystyle\sum_{g \in G}\beta_{g}\overline{g}$ two elements  of $R^{\gamma}G$, for each $k$, $0\leq k< m$, we define the \textit{generalized k-Galois form} on $R^{\gamma}G$ as 

\vr

\begin{center}
$[\alpha,\beta]_{k}=\displaystyle\sum_{g \in G}\alpha_{g}\sigma^{k}(\beta_{g})$.
\end{center}

\vr

It is not difficult to see that generalized $k$-Galois form is just the Euclidean inner product if $k=0$ and it is a sequilinear $k$-form (see \cite{ZL}). Thus, given a twisted group code $\CC$, we can define the \textit{generalized k-Galois dual code of} $\CC$ as

\begin{center}
$\CC^{\perp_{k}}=\{\beta \in R^{\gamma}G \mid [\alpha,\beta]_{k}=0, \, \forall \, \alpha \in \CC \}$.
\end{center}

Given two non-zero elements $\lambda$ and $\beta$ of $\F_{q}$, we say that a linear code $\CC$ is a $\lambda-\beta$-\textit{constacyclic} if $\CC$ is a $\lambda-$constacyclic and a $\beta$-constacyclic. Dinh, in \cite{D1}, proved if $\lambda\neq \beta$, the only $\lambda-\beta$-constacyclic codes are $\{0\}$ and $\F_{q}^{n}$.

We can generalize this result in terms of twisted group algebras. 

\begin{lem} \label{L2} {\rm{\cite [Proposition ]{D1}}}
Let $R$ be a finite commutative ring and let $C_{n}=\left\langle g \mid g^{n}=1\right\rangle$ be a cyclic group of order $n$ and $\lambda, \, \beta$ elements of ${\mathcal{U}}(R)$. Consider the twisted group algebras $R^{\gamma_{\lambda}}C_{n}$ and $R^{\gamma_{\beta}}C_{n}$ where

\vr

\begin{center}
		$\gamma_{\lambda}(g^{j},g^{k})= \left\{
		\begin{array}{lll}
			\lambda,    & {\rm{if}} \,\, j+k\geq n\\
			1, &  {\rm{if}} \,\, j+k < n.
		\end{array}
		\right.$
		$\gamma_{\beta}(g^{j},g^{k})= \left\{
		\begin{array}{lll}
			\beta,    & {\rm{if}} \,\, j+k\geq n\\
			1, &  {\rm{if}} \,\, j+k < n.
		\end{array}
		\right.$
	\end{center}
	
\vr	
	
If $\CC$ is a non-zero $\lambda$-constacyclic and also $\beta$-constacyclic code, then $\lambda=\beta$.	
\end{lem}

\begin{proof}
Let $c=\displaystyle\sum_{i=0}^{n-1}c_{i}\overline{g}^{i}$ be a non-zero element of $\CC$. Since, by hypothesis, $\CC$ is $\lambda$-constacyclic and also $\beta$-constacyclic code, we have that

\vr

$\overline{g}\cdot c= c_{0}\overline{g}\cdot\overline{1}+c_{1}\overline{g}\cdot \overline{g}+\cdots+c_{n-1}\overline{g}\cdot\overline{g^{n-1}}$

$\hspace{0.7cm}=c_{0}\cdot \overline{g}+c_{1}\overline{g^{2}}+\cdots+c_{n-1}\lambda\cdot \overline{1}$

$\hspace{0.7cm}=c_{0}\cdot \overline{g}+c_{1}\overline{g^{2}}+\cdots+c_{n-1}\beta\cdot \overline{1}$

so, $\lambda=\beta$ since the set $\{\overline{g}, g \in C_{n}\}$ is a basis of $R^{\gamma_{\lambda}}C_{n}$ and $R^{\gamma_{\beta}}C_{n}$.
\end{proof}

\begin{prop} \label{P1} 
Let $R$ be a finite commutative ring with $p^{m}$ elements, $C_{n}=\left\langle g, \, g^{n}=1\right\rangle$ be a cyclic group of order n and $R^{\gamma_{\lambda}}C_{n}$ the twisted group ring of $C_{n}$ over $R$ where

\begin{center}
		$\gamma_{\lambda}(g^{j},g^{k})= \left\{
		\begin{array}{lll}
			\lambda,    & {\rm{if}} \,\, j+k\geq n\\
			1, &  {\rm{if}} \,\, j+k < n.
		\end{array}
		\right.$
\end{center}		

Then, if $\CC$ is a $\lambda$-constacyclic code, its k-Galois dual $\CC^{\perp_{k}}$ is a $\sigma^{m-k}(\lambda^{-1})$-constacyclic code.
\end{prop}

\begin{proof}
Let $R^{\gamma_{\sigma^{k-m}(\lambda^{-1})}}C_{n}$ be the twisted group algebra with twisting defined by

\begin{center}
		$\gamma_{\sigma^{m-k}(\lambda^{-1})}(g^{j},g^{k})= \left\{
		\begin{array}{lll}
			\sigma^{m-k}(\lambda^{-1}),    & {\rm{if}} \,\, j+k\geq n\\
			1, &  {\rm{if}} \,\, j+k < n.
		\end{array}
		\right.$
\end{center}

Let $\CC$ be a $\lambda$-constacyclic code and, given any element $c=\displaystyle\sum_{i=0}^{n-1}c_{i}\overline{g}^{i}$ of the code $\CC$ and $x=\displaystyle\sum_{i=0}^{n-1}x_{i}\overline{g}^{i}$ an element of $\CC^{\perp_{k}}$ , we have
\vr

$[c,\overline{g}x]_{k}=c_{0}\sigma^{k}(x_{n-1}\sigma^{m-k}(\lambda^{-1}))+c_{1}\sigma^{k}(x_{0})+c_{2}\sigma^{k}(x_{1})+\cdots +c_{n-1}\sigma^{k}(x_{n-2})$

$\hspace{1.15cm}=c_{0}\lambda^{-1}\sigma^{k}(x_{n-1})+c_{1}\sigma^{k}(x_{0})+c_{2}\sigma^{k}(x_{1})+\cdots +c_{n-1}\sigma^{k}(x_{n-2})$


$\hspace{1.15cm}=[\overline{g}^{-1}c,x]_{k}=0$.

\vr

Then, $\overline{g}x \in \CC^{\perp_{k}}$ and the proof is completed.
\end{proof}

\begin{definition}
Let $\CC$ be a twisted group code over a finite commutative ring R. We say that $\CC$ is a linear complementary k-Galois dual code (k-Galois LCD code for shorty) if $\CC \cap \CC^{\perp_{k}}=\{0\}$.
\end{definition}

By Lemma \ref{L2} and Proposition \ref{P1}, we get

\begin{cor} {\rm{\cite [Corollary 3.3]{LFL}}} \label{OT}
If $\sigma^{m-k}(\lambda^{-1})\neq \lambda$, then any $\lambda$-constacyclic code $\CC$ over R is a k-Galois LCD code.
\end{cor}

Now, let $R$ be a finite commutative ring. We say $R$ is a {\textit{chain ring}} if the set of all ideals of $R$ is a chain under set-theoretic inclusion. Futhermore, the following conditions are equivalent

(i) $R$ is a local ring and the maximal ideal $M$ of $R$ is principal,

(ii) $R$ is a local principal ideal ring,

(iii) $R$ is a chain ring.

Finally, if $R$ is a finite commutative chain ring, we have that for any $k$-Galois LCD twisted group code, $|\CC|\cdot |\CC^{\perp_{k}}|=|R|^{n}$ and $\CC$ is free with $\CC \oplus \CC^{\perp_{k}}=R^{\gamma}G$ which implies $rank_{R}(\CC) + rank_{R}(\CC^{\perp_{k}})=|G|=n$. See, for example, \cite{BBBFM} Theorem 2.


Let $R$ be a commutative ring with identity and let $G$ be a group. Consider the following mapping $^{\ast} : R^{\gamma}G \longrightarrow R^{\gamma}G $ given by $\left(\displaystyle\sum_{g \in G} \alpha_{g}\overline{g}\right)^{\ast}=\displaystyle\sum_{g \in G} \alpha_{g}\overline{g}^{-1}$.

\begin{definition}
Let R be a commutative ring with identity and let G be a group. The mapping $^{\ast} : R^{\gamma}G \longrightarrow R^{\gamma}G $ given by $\left(\displaystyle\sum_{g \in G} \alpha_{g}\overline{g}\right)^{\ast}=\displaystyle\sum_{g \in G} \alpha_{g}\overline{g}^{-1}$ with $\gamma(g,h)^{2}=1$, is called the classical involution of $R^{\gamma}G$.
\end{definition}

\begin{lem} \label{L1}
Let $R$ be a finite commutative chain ring with characteristc $p$, $C_{n}=\left\langle g, \, g^{n}=1\right\rangle$ be a cyclic group of order n and $R^{\gamma_{\lambda}}C_{n}$ the twisted group algebra of $C_{n}$ over $R$ where

\begin{center}
		$\gamma_{\lambda}(g^{j},g^{k})= \left\{
		\begin{array}{lll}
			\lambda,    & {\rm{if}} \,\, j+k\geq n\\
			1, &  {\rm{if}} \,\, j+k < n.
		\end{array}
		\right.$
\end{center}		

Given two arbitrary elements $\alpha=\displaystyle\sum_{i=0}^{n-1}\alpha_{i}\overline{g}^{i}$ and $\beta=\displaystyle\sum_{i=0}^{n-1}\beta_{i}\overline{g}^{i}$ of $R^{\gamma_{\lambda}}C_{n}$, let us denote by $\sigma^{k}(\beta)$ the element $\displaystyle\sum_{i=0}^{n-1}\sigma^{k}(\beta_{i})\overline{g}^{i}$. If $\alpha\left(\sigma^{k}(\beta)\right)^{\ast}=0$ and $\lambda^{2}=1$, then $[\alpha,\beta]_{k}=0$.
\end{lem}

\begin{proof}
It is not difficult to see, the coefficient of $1=\overline{1}$ in the product $\alpha\left(\sigma^{k}(\beta)\right)^{\ast}$ is exactly $[\alpha,\beta]_{k}$. Since, by hypothesis, $\alpha\left(\sigma^{k}(\beta)\right)^{\ast}=0$, we have $[\alpha,\beta]_{k}=0$. 
\end{proof}

\section{Euclidean Constacyclic LCD codes.}

In this section, we shall characterize constacyclic LCD codes in terms of its idempotent generator with respect to Euclidean inner product.

Let $R$ be a finite commutative chain ring with $p^{m}$ elements, $C_{n}=\left\langle g, \, g^{n}=1\right\rangle$ be a cyclic group of order $n$ and $R^{\gamma_{\lambda}}C_{n}$ the twisted group algebra of $C_{n}$ over $R$ where

\begin{center}
		$\gamma_{\lambda}(g^{j},g^{k})= \left\{
		\begin{array}{lll}
			\lambda,    & {\rm{if}} \,\, j+k\geq n\\
			1, &  {\rm{if}} \,\, j+k < n.
		\end{array}
		\right.$
\end{center}		

\noindent for some unit of R. Given $\alpha=\displaystyle\sum_{g \in C_{n}}\alpha_{g}\overline{g}$, $\beta=\displaystyle\sum_{g \in C_{n}}\beta_{g}\overline{g}$ two elements  of $R^{\gamma_{\lambda}}C_{n}$, we define the \textit{Euclidean inner product} on $R^{\gamma_{\lambda}}C_{n}$ as 

\vr

\begin{center}
$[\alpha,\beta]=\displaystyle\sum_{g \in G}\alpha_{g}\beta_{g}$.
\end{center}

Let $\CC$ be a constacyclic code over $R^{\gamma}C_{n}$, that is, an ideal of $R^{\gamma_{\lambda}}C_{n}$ . It is well-known that the set $\CC^{\perp}=\{x \in R^{\gamma}C_{n} \mid [x,\alpha]=0, \forall \alpha \in \CC\}$ is an ideal in the twisted group algebra $R^{\gamma_{\lambda^{-1}}}C_{n}$ where

\begin{center}
		$\gamma_{\lambda^{-1}}(g^{j},g^{k})= \left\{
		\begin{array}{lll}
			\lambda^{-1},    & {\rm{if}} \,\, j+k\geq n\\
			1, &  {\rm{if}} \,\, j+k < n.
		\end{array}
		\right.$
\end{center}		

\begin{definition}
Let $\CC$ be a constacyclic code over a finite commutative ring $R$. We say that $\CC$ is a linear complementary dual code (LCD code for shorty) if $\CC \cap \CC^{\perp}=\{0\}$.
\end{definition}

Notice that, if $\lambda^{2}\neq 1$, the Corollary \ref{OT} shows us any $\lambda$-constacyclic code is a LCD code. 

\begin{prop} \label{LCDE}
Let $R$ be a finite commutative chain ring with $p^{m}$ elements, with p rational prime, $C_{n}=\left\langle g, \, g^{n}=1\right\rangle$ be a cyclic group of order n and $R^{\gamma_{\lambda}}C_{n}$ the twisted group algebra of $C_{n}$ over $R$ where

\begin{center}
		$\gamma_{\lambda}(g^{j},g^{k})= \left\{
		\begin{array}{lll}
			\lambda,    & {\rm{if}} \,\, j+k\geq n\\
			1, &  {\rm{if}} \,\, j+k < n.
		\end{array}
		\right.$
\end{center}		

\noindent for some unit of R. If $\lambda^{2}=1$, then $\CC$ is a $\lambda$-constacyclic LCD code if, and only if, $\CC$ is generated by an idempontent e such that $e=e^{\ast}$.
\end{prop}

\begin{proof}
First of all, if $\CC$ is a $\lambda$-constacyclic code which is also LCD, we have the following decomposition of ideals $R^{\gamma_{\lambda}}C_{n}=\CC \oplus \CC^{\perp}$ since $\lambda^{2}=1$. So , it is well-know, there exist idempotents $e$ and $f$ such that $1=e+f$, $e\cdot f=0$, $\CC=\left\langle e\right\rangle$ and $\CC^{\perp}=\left\langle f\right\rangle$. 

Since $[e,1-e]=0$, we have that $[1,e^{\ast}(1-e)]=0$. Now, given $a=\displaystyle\sum_{i=0}^{n-1}a_{i}\overline{g}^{i}$ and $b=\displaystyle\sum_{i=0}^{n-1}b_{i}\overline{g}^{i}$ two arbitrary elements of $R^{\gamma_{\lambda}}C_{n}$, then

\begin{center}
$[\overline{g}a,\overline{g}b]=a_{0}b_{0}+a_{1}b_{1}+\cdots+a_{n-2}b_{n-2}+(a_{n-1}\lambda)(b_{n-1}\lambda)=[a,b]$
\end{center} 
 
So $[\overline{g},\overline{g}e^{\ast}(1-e)]=0$, for all $g \in G$. Since the Euclidean inner product is non-degenerated, we get that $e^{\ast}(1-e)=0$ and $e^{\ast}=e^{\ast}e$, which implies that $e=(e^{\ast})^{\ast}=(e^{\ast}e)^{\ast}=e^{\ast}e=e^{\ast}$. 

On the other hand, if $e$ is an idempotent such that $e=e^{\ast}$ and $\CC=\left\langle e\right\rangle$, then, $e(1-e)^{\ast}=e(1-e^{\ast})=e(1-e)=0$ and, by Lemma \ref{L1}, $[e,1-e]=0$. Writing $1=e+(1-e)$ we have $R^{\gamma_{\lambda}}C_{n}=\CC \oplus R^{\gamma_{\lambda}}C_{n}(1-e)$. Since $[e,1-e]=0$ and  $rank_{R}(\CC) + rank_{R}(R^{\gamma_{\lambda}}C_{n}(1-e))=n$, we conclude that $\CC^{\perp}=R^{\gamma_{\lambda}}C_{n}(1-e)$.

Thus, $\CC$ is a $\lambda$-constacyclic LCD code. 
\end{proof}

Thus, we obtain, as a corollary, the Theorem 3.1 in \cite{CW}.

\begin{cor} \label{C1} 
Let $R$ be a finite commutative chain ring with $p^{m}$ elements, $C_{n}=\left\langle g, \, g^{n}=1\right\rangle$ be a cyclic group of order n. A cyclic code $\CC$ is a LCD code with respect to Euclidean inner product, if, and only if,  $\CC=\left\langle e \right\rangle$ such that $e^{2}=e$ and $e=e^{\ast}$.
\end{cor}

\section{Hermitian LCD constacyclic codes over $\Z_{p^{2m}}$}

In this section, we shall characterize constacyclic LCD codes over the finite chain ring $\Z_{p^{2m}}$  in terms of its idempotent generator with respect to Hermitian inner product. In this case, the automorphism group of $\Z_{p^{2m}}$ is generated by the Frobenius automorphism.

Let $R=\Z_{p^{2m}}$ be a finite commutative chain ring with $p^{2m}$ elements, $C_{n}=\left\langle g, \, g^{n}=1\right\rangle$ be a cyclic group of order $n$ and $R^{\gamma_{\lambda}}C_{n}$ the twisted group algebra of $C_{n}$ over $R$ where

\begin{center}
		$\gamma_{\lambda}(g^{j},g^{k})= \left\{
		\begin{array}{lll}
			\lambda,    & {\rm{if}} \,\, j+k\geq n\\
			1, &  {\rm{if}} \,\, j+k < n.
		\end{array}
		\right.$
\end{center}		

\noindent for some unit of R. Given $\alpha=\displaystyle\sum_{g \in C_{n}}\alpha_{g}\overline{g}$, $\beta=\displaystyle\sum_{g \in C_{n}}\beta_{g}\overline{g}$ two elements  of $R^{\gamma_{\lambda}}C_{n}$, we define the \textit{Hermitian inner product} on $R^{\gamma_{\lambda}}C_{n}$ by 

\vr

\begin{center}
$[\alpha,\beta]_{m}=\displaystyle\sum_{g \in G}\alpha_{g}\beta^{p^{m}}_{g}$.
\end{center}

Let $\CC$ be a constacyclic code over $R^{\gamma}C_{n}$, that is, an ideal of $R^{\gamma_{\lambda}}C_{n}$ . It is well-known that the set $\CC^{\perp}=\{x \in R^{\gamma}C_{n} \mid [x,\alpha]_{m}=0, \forall \alpha \in \CC\}$ is an ideal in the twisted group algebra $R^{\gamma_{\lambda^{-p^{m}}}}C_{n}$ where

\begin{center}
		$\gamma_{\lambda^{-p^{m}}}(g^{j},g^{k})= \left\{
		\begin{array}{lll}
			\lambda^{-p^{m}},    & {\rm{if}} \,\, j+k\geq n\\
			1, &  {\rm{if}} \,\, j+k < n.
		\end{array}
		\right.$
\end{center}		

\begin{definition}
Let $\CC$ be a constacyclic code over a finite commutative ring $R$. We say that $\CC$ is a linear complementary dual code (LCD code for shorty) if $\CC \cap \CC^{\perp_{m}}=\{0\}$.
\end{definition}

Notice that, if $\lambda^{p^{m}+1}\neq 1$, the Corollary \ref{OT} shows us any $\lambda$-constacyclic code is a LCD code. Also, if we take $\lambda$ with $\lambda^{2}=1$, $p$ has to be an odd prime. If $p=2$, both conditions $\lambda^{2}=1$ and $\lambda^{p^{m}+1}=1$ imply $\lambda=1$. So, in this case the 2-cocycle is trivial then $R^{\gamma_{\lambda}}C_{n}$ is the ordinary group algebra.

\begin{thrm} \label{T1} 
Let $R$ be a finite commutative chain ring with $p^{2m}$, $C_{n}=\left\langle g, \, g^{n}=1\right\rangle$ be a cyclic group of order n and $R^{\gamma_{\lambda}}C_{n}$ the twisted group algebra of $C_{n}$ over $R$ where

\begin{center}
		$\gamma_{\lambda}(g^{j},g^{k})= \left\{
		\begin{array}{lll}
			\lambda,    & {\rm{if}} \,\, j+k\geq n\\
			1, &  {\rm{if}} \,\, j+k < n.
		\end{array}
		\right.$
\end{center}		

Let e be an idempotent of $R^{\gamma_{\lambda}}C_{n}$ and $\lambda^{2}=1$. Then $e=e(e^{(p^{m})})^{\ast}$ if, and only if $[e,1-e]_{m}=0$. 
\end{thrm}

\begin{proof}
Suppose that $e$ is an idempotent such that $e=e(e^{(p^{m})})^{\ast}$. Then,\\ $0=e-e(e^{(p^{m})})^{\ast}=e(1-(e^{(p^{m})})^{\ast})=e\left((1-e)^{(p^{m})}\right)^{\ast}$ and, by Lemma \ref{L1}, \\ $[e,1-e]_{m}=0$.

On the other hand, if $[e,1-e]_{m}=0$, we have that $[1,e^{\ast}(1-e)^{(p^{m})}]=0$. Now, given $a=\displaystyle\sum_{i=0}^{n-1}a_{i}\overline{g}^{i}$ and $b=\displaystyle\sum_{i=0}^{n-1}b_{i}\overline{g}^{i}$ two arbitrary elements of $R^{\gamma_{\lambda}}C_{n}$, since $\lambda^{p^{m}}=\lambda^{-1}$, then 

\vr

$[\overline{g}a,\overline{g}b]_{m}=a_{0}b_{0}^{p^{m}}+a_{1}b_{1}^{p^{m}}+\cdots+a_{n-2}b_{n-2}^{p^{m}}+(a_{n-1}\lambda)(b_{n-1}^{p^{m}}\lambda^{p^{m}})=[a,b]_{}$

\vr
 
So $[\overline{g},\overline{g}e^{\ast}(1-e)^{(p^{m})}]=0$, for all $g \in G$. Since the Hermitian form is non-degenerated, we get that $e^{\ast}(1-e)^{(p^{m})}=0$ so $e^{\ast}=e^{\ast}e^{(p^{m})}$. Then, $e=(e^{\ast})^{\ast}=(e^{\ast}e^{(p^{m})})^{\ast}=e(e^{(p^{m})})^{\ast}$. 
\end{proof}

Now, we have the following

\begin{prop} \label{LCDH}
Let $R$ be a finite commutative chain ring with identity and $p^{2m}$, with p a rational prime, $C_{n}=\left\langle g, \, g^{n}=1\right\rangle$ be a cyclic group of order n and $R^{\gamma_{\lambda}}C_{n}$ the twisted group algebra of $C_{n}$ over $R$ where

\begin{center}
		$\gamma_{\lambda}(g^{j},g^{k})= \left\{
		\begin{array}{lll}
			\lambda,    & {\rm{if}} \,\, j+k\geq n\\
			1, &  {\rm{if}} \,\, j+k < n.
		\end{array}
		\right.$
\end{center}		

If $\lambda^{2}=1$, then $\CC$ is a $\lambda$-constacyclic code  generated by an idempontent e such that $e=e(e^{(p^{m})})^{\ast}$ if, and only if $\CC$ is Hermitian LCD code.
\end{prop}

\begin{proof}
First of all, if $\CC$ is a $\lambda$-constacyclic code which is also $k$-Galois LCD, we have the following decomposition of ideals $R^{\gamma}C_{n}=\CC \oplus \CC^{\perp_{k}}$ since $\lambda^{1+p^{m}}=1$.

It is well-know that there exist idempotents $e$ and $f$ such that $1=e+f$, $e\cdot f=0$, $\CC=\left\langle e\right\rangle$ and $\CC^{\perp_{m}}=\left\langle f\right\rangle$. Then, writing $f=1-e$, we get $[e,1-e]_{m}=0$, so, by Theorem \ref{T1}, the equality $e=e(e^{(p^{m})})^{\ast}$ holds. 

If $\CC$ is generated by an idempontent $e$ such that $e=e(e^{(p^{m})})^{\ast}$, by Theorem \ref{T1}, we have $[e,1-e]_{m}=0$. Then, writing $1=e+(1-e)$ we have $R^{\gamma}C_{n}=\CC \oplus R^{\gamma}C_{n}(1-e)$. Since $[e,1-e]_{m}=0$ and  $rank_{R}(\CC) + rank_{R}(R^{\gamma}C_{n}(1-e))=n$, we conclude that $\CC^{\perp_{m}}=R^{\gamma}C_{n}(1-e)$. Thus, $\CC$ is a $\lambda$-constacyclic Hermitian LCD code. 
\end{proof}

\begin{prop}
Let $R$ be a finite commutative chain ring with identity and $p^{2m}$ elements, with p a rational prime, $C_{n}=\left\langle g, \, g^{n}=1\right\rangle$ be a cyclic group of order n and $R^{\gamma_{\lambda}}C_{n}$ the twisted group algebra of $C_{n}$ over $R$ where

\begin{center}
		$\gamma_{\lambda}(g^{j},g^{k})= \left\{
		\begin{array}{lll}
			\lambda,    & {\rm{if}} \,\, j+k\geq n\\
			1, &  {\rm{if}} \,\, j+k < n.
		\end{array}
		\right.$
\end{center}		

Suppose that $\lambda^{2}=1$ and  $\CC$ is a $\lambda$-constacyclic code  generated by an idempontent e. Then, $C \subset C^{\perp_{m}}$ if, and only if, $e(e^{(p^{m})})^{\ast}=0$.
\end{prop}

\begin{proof}
Let $\CC$ be $\lambda$-constacyclic code such that $\lambda^{2}=1$ and $\CC=\left\langle e \right\rangle$ with $e^{2}=e$. If $C \subset C^{\perp_{m}}$, then $0=[e,e]_{m}=[1,e^{\ast}(e^{(p^{m})})]_{m}$. Thus, $[\overline{g},\overline{g}e^{\ast}(e^{(p^{m})})]_{m}=0$, for all $g \in C_{n}$. It implies $e^{\ast}(e^{(p^{m})})=0$ and $e(e^{(p^{m})})^{\ast}=(e^{\ast}(e^{(p^{m})}))^{\ast}=0$.

On the other hand, if $e(e^{(p^{m})})^{\ast}=0$, then $[\overline{g}e, \overline{h}e]_{m}=[\overline{g},\overline{h}e^{\ast}(e^{(p^{m})})]_{m}=0$, for all $\overline{g}$, $\overline{h} \in C_{n}$. So, $[\alpha,\beta]_{m}=0$, for all $\alpha$ and $\beta \in \CC$ and the proof is completed.
\end{proof}

\vr

\noindent \textit{Example:} Let $R=\Z_{9}$ be the ring of integers modulo 9, $C_{11}$ be the cyclic group of order 11 and $R^{\gamma}C_{11}$ be the twisted group ring where

\begin{center}
		$\gamma(g^{j},g^{k})= \left\{
		\begin{array}{lll}
			8,    & {\rm{if}} \,\, j+k\geq 11\\
			1, &  {\rm{if}} \,\, j+k < 11.
		\end{array}
		\right.$
\end{center}		

So, in this case, we have $\overline{g}^{11}=8$. Taking 

$e= 4\overline{g}^{10} + 5\overline{g}^9 + 4\overline{g}^8 + 5\overline{g}^7 + 4\overline{g}^6 + 5\overline{g}^5 + 4\overline{g}^4 + 5\overline{g}^3 + 4\overline{g}^2 + 5\overline{g} + 5$, \\ we have $e^{2}=e$,  

$e^{(3)}=\overline{g}^{10} + 8\overline{g}^9 + \overline{g}^8 + 8\overline{g}^7 + \overline{g}^6 + 8\overline{g}^5 + \overline{g}^4 + 8\overline{g}^3 + \overline{g}^2 + 8\overline{g} + 8$ and 

$(e^{(3)})^{\ast}=8\overline{g}^{10} + \overline{g}^9 + 8\overline{g}^8 + \overline{g}^7 + 8\overline{g}^6 + \overline{g}^5 + 8\overline{g}^4 + \overline{g}^3 + 8\overline{g}^2 + \overline{g} + 8$.

Finally, it is not difficulty to see, $e(e^{(3)})^\ast=0$, so we have found a self-orthogonal code with respect to Hermitian inner product.


\section{Some good LCD codes}

We shall recall the following results.

\begin{prop} \label{BA1} {\rm{\cite [Proposition 7.20]{BA2}}}
Let M be a free module over a commutative ring R with base of n elements. Then (1) any base has cardinality n, (2) any set of generators contains at least n elements, (3) any set of n generators is a base.
\end{prop}

It is well known the number of any bases of a free module $M$ over a commutative ring is called \textit{rank of M}.

\begin{prop} \label{BA2} {\rm{\cite [Corollary page 415]{BA2}}}
Let M be a finitely generated R-module. The $x_{1}, x_{2}, \cdots, x_{m}$ generate M if and only if the cosets $\overline{x_{1}}=x_{1}+(J(R))M, \overline{x_{2}}=x_{2}+(J(R))M, \cdots, \overline{x_{m}}=x_{m}+(J(R))M$ generate $M/(J(R))M$ as $R/J(R)$-module.
\end{prop}

We say that a $R$-module $M$ is \textit{projective} if $M$ is a direct sommand of a free $R$-module. Thus, 

\begin{prop} \label{BA3} {\rm{\cite [Theorem 7.7]{BA2}}}
Any finitely generated projective module over a commutative local ring is free.
\end{prop}

\begin{thrm}[Generalized Singleton Bound] \label{MDR1} {\rm{\cite [Theorem 4.12]{DO}}}
Let $\CC$ be a linear code of length n over a finite commutative chain ring R. Let k = min\{$\ell$ $\mid$ there exists a monommorphism from $\CC$ to $R^{\ell}$ as R modules \}. Then

\begin{center}
$d_{H}(\CC) \leq n - k + 1.$
\end{center}
\end{thrm}

If a linear code $\CC$ meets the bound above, we say the code is a {\textit{Maximum Distance}} with respect the {\textit{Rank} (MDR)} code. Also, when $|\CC|=q^{k}$, we have $d_{H}(\CC) \leq n - k + 1$ which is called {\textit{Singleton Bound}} and if $\CC$ meets the Singleton Bound, we say $\CC$ is a {\textit{Maximum Distance Separable} (MDS)} code.

Notice that, an MDR code is not necessarily MDS, see, for example, \cite{DO} page 57. But, if $R$ is a finite commutative chain ring, we have

\begin{cor} \label{MDR2} {\rm{\cite [Corollary 4.4]{DO}}}
Let R be a finite commutative chain ring and let $\CC$ be a linear code over R. The code $\CC$ is an MDS code if and only if $\CC$ is an MDR code and $\CC$ is free. 
\end{cor}

Now, let $R$ be a finite commutative chain ring with $R/J(R) = \F_{p^{s}}$ and let $c=\displaystyle\sum_{i=0}^{n-1}c_{i}\overline{g}^{i}$ be an element of  $R^{\gamma}C_{n}$ . Consider the following mapping 

\begin{center}
$\vartheta: R^{\gamma}C_{n} \longrightarrow \F_{p^{s}}^{\overline{\gamma}}C_{n}$ 
\end{center}

given by $\vartheta(c)=\displaystyle\sum_{i=0}^{n-1}\left(c_{i}+J(R)\right)\overline{g}^{i}$ where $\overline{\gamma}=\gamma+J(R)$. 

The mapping $\vartheta$ above defined has the following properties:

\begin{enumerate}
	\item $\vartheta(xy)=\vartheta(x)\vartheta(y)$;
	
	\item $w_{H}(x)\geq w_{H}(\vartheta(x))$, for all $x,y \in R^{\gamma}C_{n}$.
\end{enumerate}

Notice that, $(a\overline{g})(b\overline{h})=ab\gamma(g,h)\overline{gh}$, then $\vartheta((a\overline{g})(b\overline{h}))=((ab\gamma(g,h))+J(R))\overline{gh}=(a+J(R))(b+J(R))(\gamma(g,h)+J(R))\overline{gh}=((a+J(R))\overline{g})((b+J(R))\overline{h})=\vartheta(a\overline{g})\vartheta(b\overline{h})$. This proves the first property.

The second property is trivial.

\vr

The weights of the codes presented in this section were computed using SageMath. 

\noindent {\textit{Example 1:}} Let $R=\Z_{4}$ be the ring of integers modulo 4, $C_{19}$ be the cyclic group of order 19 and $R^{\gamma}C_{19}$ be the twisted group ring where

\begin{center}
		$\gamma(g^{j},g^{k})= \left\{
		\begin{array}{lll}
			3,    & {\rm{if}} \,\, j+k\geq 19\\
			1, &  {\rm{if}} \,\, j+k < 19.
		\end{array}
		\right.$
\end{center}		

So, in this case, we have $\overline{g}^{19}=3$. Now, taking the element 

{\small{
$e=\overline{g}^{18} + 3\overline{g}^{17} + \overline{g}^{16} + 3\overline{g}^{15} + \overline{g}^{14} + 3\overline{g}^{13} + \overline{g}^{12} + 3\overline{g}^{11} + \overline{g}^{10} + 3\overline{g}^9 + \overline{g}^8 + 3\overline{g}^7 + \overline{g}^6 + 3\overline{g}^5$ 

\hspace{0.25cm} $ + \overline{g}^4 + 3\overline{g}^3 + \overline{g}^2 + 3\overline{g} + 2$ }}, we have $e=e^{2}=e^{\ast}$ and, by Proposition \ref{LCDE}, the code generated by $e$ is a LCD code of rank 18 and weight 2. Consequently, we have found a MDR code which has the same weight of the best [19,18]-code over the field $\F_{4}$. Notice also this code is MDS, by Corollary \ref{MDR2}.

Finally, we can consider the idempotent $f=1-e$ such that $f^{\ast}=f$ and $\CC^{\perp}=\left\langle f\right\rangle$ which implies that code generated by $f$ is also LCD of rank 1 and weight 19. So, we have found a MDR code which has the same weight of the best [19,1]-code over the field $\F_{4}$. Notice also this code is MDS, by Corollary \ref{MDR2}.

\vr

\noindent {\textit{Example 2:}} Let $R=\Z_{8}$ be the ring of integers modulo 8, $C_{3}$ be the cyclic group of order 3 and $R^{\gamma}C_{3}$ be the twisted group ring where

\begin{center}
		$\gamma(g^{j},g^{k})= \left\{
		\begin{array}{lll}
			3,    & {\rm{if}} \,\, j+k\geq 3\\
			1, &  {\rm{if}} \,\, j+k < 3.
		\end{array}
		\right.$
\end{center}		

So, in this case, we have $\overline{g}^{3}=3$. Now, taking the element $e=3\overline{g}^{2}+\overline{g}+3$, we have we have $e=e^{2}=e^{\ast}$.

	
	
	
	
  
Thus, $e$ is an idempotent such that $e^{\ast}=e$. By Proposition \ref{LCDE}, the code generated by $e$ is a LCD code of rank 1 and weight 3. Consequently, we have found a MDR code which has the same weight of the best [3,1]-code over the field $\F_{8}$. Notice also this code is MDS, by Corollary \ref{MDR2}.	

Finally, we can consider the idempotent $f=1-e=5\overline{g}^{2}+7\overline{g}+6$ such that $f^{\ast}=f$ and $\CC^{\perp}=\left\langle f\right\rangle$ which implies that code generated by $f$ is also LCD of rank 2 and weight 2. So, we have found a MDR code which has the same weight of the best [3,2]-code over the field $\F_{8}$. Notice also this code is MDS, by Corollary \ref{MDR2}.	

Also, those code generated by $e$ and $f$ are codes of {\textit{constant weight}}, that is, all non-zero codewords have the same weight.

\vr





\noindent \textit{Example 3:} Let $C_{21}=\left\langle g, \,\, g^{21}=1\right\rangle$ be a cyclic group of order 21 and let $R=\F_{5}$ be a finite field with 5 elements. Consider the twisted group algebra $\F_{5}^{\gamma_{4}}C_{21}$ where

\begin{center}
		$\gamma(g^{j},g^{k})= \left\{
		\begin{array}{lll}
			4,    & {\rm{if}} \,\, j+k\geq 21\\
			1, &  {\rm{if}} \,\, j+k < 21.
		\end{array}
		\right.$
\end{center}

Thus in this case, we have $\overline{g}^{21}=4$. Finally, taking the element 

\vr

\noindent $e=\overline{g}^{20} + 4\overline{g}^{18} + 4\overline{g}^{17} + \overline{g}^{16} + \overline{g}^{15} +2\overline{g}^{14}+ 4\overline{g}^{12} + \overline{g}^9  + 3\overline{g}^7 + 4\overline{g}^6 + 4\overline{g}^5 + \overline{g}^4 + \overline{g}^3 + 4\overline{g} +1 $, it is not difficult to see $e^{2}=e=e^{\ast}$. 









So, by Proposition \ref{LCDE}, the code $\CC$ generated by $e$ is a LCD code of dimension 6 and weight 12 which are exactly the parameters of the best [21,6] code known.

\vr

\noindent \textit{Example 4:} Let $C_{10}=\left\langle g, \,\, g^{10}=1\right\rangle$ be a cyclic group of order 10 and let $\F_{3}$ be a finite field with 3 elements. Consider the twisted group algebra $\F_{3}^{\gamma_{2}}C_{10}$ where

\begin{center}
		$\gamma(g^{j},g^{k})= \left\{
		\begin{array}{lll}
			2,    & {\rm{if}} \,\, j+k\geq 10\\
			1, &  {\rm{if}} \,\, j+k < 10.
		\end{array}
		\right.$
\end{center}

Thus in this case, we have $\overline{g}^{10}=2$. Finally, taking the elements $e=\overline{g}^{8}+2\overline{g}^{6}+\overline{g}^{4}+2\overline{g}^2+2$ and $f=2\overline{g}^8 + \overline{g}^6 + 2\overline{g}^4 + \overline{g}^2 + 2 $.

It is not difficult to see $e^{2}=e=e^{\ast}$ 


So, by Proposition \ref{LCDE}, the code $\CC$ generated by $e$ is a LCD code of dimension of dimension 8 and weight 2 which are exactly the parameters of the best [10,8] code known.

Finally, notice that $f=1-e$, so it is also an idempotent with $f^{\ast}=f$ and the code generated by $f$ is LCD of dimension 2 and weight 5 and the best [10,2] code known has weight 7.

\vr

\noindent \textit{Example 5:} Let $R=\Z_{16}$ be the ring of integers modulo 16, $C_{33}$ be the cyclic group of order 33 and $R^{\gamma}C_{33}$ be the twisted group ring where

\begin{center}
		$\gamma(g^{j},g^{k})= \left\{
		\begin{array}{lll}
			7,    & {\rm{if}} \,\, j+k\geq 33\\
			1, &  {\rm{if}} \,\, j+k < 33.
		\end{array}
		\right.$
\end{center}		

So, in this case, we have $\overline{g}^{33}=7$. Taking

$e=9\overline{g}^{32} + 15\overline{g}^{31} + 2\overline{g}^{30} + 15\overline{g}^{29} + 8\overline{g}^{28} + 14\overline{g}^{27} + 8\overline{g}^{26} + 15\overline{g}^{25} + 2\overline{g}^{24} + 8\overline{g}^{23} + 11\overline{g}^{22} + 14\overline{g}^{21} + 8\overline{g}^{20} + 8\overline{g}^{19} + 2\overline{g}^{18} + 15\overline{g}^{17} + 9\overline{g}^{16} + 14\overline{g}^{15} + 8\overline{g}^{14} + 8\overline{g}^{13} + 2\overline{g}^{12} + 13\overline{g}^{11} + 8\overline{g}^{10} + 14\overline{g}^9 + 9\overline{g}^8 + 8\overline{g}^7 + 2\overline{g}^6 + 8\overline{g}^5 + 9\overline{g}^4 + 14\overline{g}^3 + 9\overline{g}^2 + 15\overline{g} + 13$, we have $e=e^{2}=e^{\ast}$ and, by Proposition \ref{LCDE}, the code $\CC$ generated by $e$ is a LCD code of dimension of dimension 13.

Notice that, the element ${\tilde{e}}=\vartheta(e)= g^{32} + g^{31} + g^{29} + g^{25} + g^{22} + g^{17} + g^{16} + g^{11} + g^8 + g^4 + g^2 + g + 1  \in \F_{2}C_{33}$ is also an idempotent with ${\tilde{e}}^{\ast}={\tilde{e}}$ and the code generated by ${\tilde{e}}$ is a LCD code of dimension 13 and weight 10 which are exactly the parameters of the best [33,13] code known over $\F_{2}$. Also, the code generated by $1-{\tilde{e}}$ is a LCD code of dimension 20 and weight 6 which are exactly the parameters of the best [33,20] code known over $\F_{2}$.




\section{k-Galois constacyclic LCD codes}

In this section, we shall prove some results about $k$-Galois constacyclic LCD codes over a commutative chain ring $R$ with characteristic $p$. Let $\sigma$ be an automorphism of $R$ of order $m$.

\begin{thrm} \label{T1}
Let $R$ be a finite commutative chain ring with characteristic p, $C_{n}=\left\langle g, \, g^{n}=1\right\rangle$ be a cyclic group of order n and $R^{\gamma_{\lambda}}C_{n}$ the twisted group algebra of $C_{n}$ over $R$ where

\begin{center}
		$\gamma_{\lambda}(g^{j},g^{k})= \left\{
		\begin{array}{lll}
			\lambda,    & {\rm{if}} \,\, j+k\geq n\\
			1, &  {\rm{if}} \,\, j+k < n.
		\end{array}
		\right.$
\end{center}		

Let e be an idempotent of $R^{\gamma_{\lambda}}C_{n}$ and $\lambda^{2}=1$. Then $e=e(\sigma^{k}(e))^{\ast}$ if, and only if $[e,1-e]_{k}=0$. 
\end{thrm}

\begin{proof}
Suppose that $e$ is an idempotent such that $e=e(\sigma^{k}(e))^{\ast}$. Then,

$0=e-e(\sigma^{k}(e))^{\ast}=e(1-(\sigma^{k}(e))^{\ast})=e\left(\sigma^{k}((1-e)\right)^{\ast}$ and, by Lemma \ref{L1}, $[e,1-e]_{k}=0$.

On the other hand, if $[e,1-e]_{k}=0$, we have that $[1,e^{\ast}\sigma^{k}(1-e)]=0$. Now, given $a=\displaystyle\sum_{i=0}^{n-1}a_{i}\overline{g}^{i}$ and $b=\displaystyle\sum_{i=0}^{n-1}b_{i}\overline{g}^{i}$ two arbitrary elements of $R^{\gamma_{\lambda}}C_{n}$, since $\sigma^{m-k}(\lambda)=\lambda$, then 

$[\overline{g}a,\overline{g}b]_{k}=a_{0}\sigma^{k}(b_{0})+a_{1}\sigma^{k}(b_{1})+\cdots+a_{n-2}\sigma^{k}(b_{n-2})+(a_{n-1}\lambda)(\sigma^{k}(b_{n-1}\lambda))$

\hspace{1.15cm} $=[a,b]_{k}$
 
So $[\overline{g},\overline{g}e^{\ast}\sigma^{k}(1-e)]=0$, for all $g \in G$. Since the $k$-Galois form is non-degenerated, we get that $e^{\ast}\sigma^{k}(1-e)=0$ so $e^{\ast}=e^{\ast}\sigma^{k}(e)$. Then, $e=(e^{\ast})^{\ast}=(e^{\ast}\sigma^{k}(e))^{\ast}=e(\sigma^{k}(e))^{\ast}$. 
\end{proof}

Notice that, if $\sigma^{m-k}(\lambda^{-1})\neq \lambda$, the Corollary \ref{OT} shows us any $\lambda$-constacyclic code is a LCD code. Also, if $\lambda^{2}=1$ and $\sigma^{m-k}(\lambda)\neq \lambda$, any $\lambda$-constacyclic code is a LCD code.   


Now, for the next result, we shall consider $\lambda^{2}=1$ and $\sigma^{m-k}(\lambda)=\lambda$.

\begin{prop} \label{GL}
Let $R$ be a finite commutative chain ring with identity and characteristic p, a rational prime, $C_{n}=\left\langle g, \, g^{n}=1\right\rangle$ be a cyclic group of order n and $R^{\gamma_{\lambda}}C_{n}$ the twisted group algebra of $C_{n}$ over $R$ where

\begin{center}
		$\gamma_{\lambda}(g^{j},g^{k})= \left\{
		\begin{array}{lll}
			\lambda,    & {\rm{if}} \,\, j+k\geq n\\
			1, &  {\rm{if}} \,\, j+k < n.
		\end{array}
		\right.$
\end{center}		

If $\lambda^{2}=1$, then $\CC$ is a $\lambda$-constacyclic code  generated by an idempontent e such that $e=e(\sigma^{k}(e))^{\ast}$ if, and only if $\CC$ is k-Galois LCD code.
\end{prop}

\begin{proof}
If $\CC$ is a $\lambda$-constacyclic code which is also $k$-Galois LCD, we have $\sigma^{m-k}(\lambda^{-1})=\sigma^{m-k}(\lambda)=\sigma^{k}(\lambda)=\lambda$ since $\lambda^{2}=1$.


It is well-know that there exist idempotents $e$ and $f$ such that $1=e+f$, $e\cdot f=0$, $\CC=\left\langle e\right\rangle$ and $\CC^{\perp_{k}}=\left\langle f\right\rangle$. Then, writing $f=1-e$, we get $[e,1-e]_{k}=0$, so, by Theorem \ref{T1}, the equality $e=e(\sigma^{k}(e))^{\ast}$ holds. 

If $\CC$ is generated by an idempontent $e$ such that $e=e(\sigma^{k}(e))^{\ast}$, by Theorem \ref{T1}, we have $[e,1-e]_{k}=0$. Then, writing $1=e+(1-e)$ we have $R^{\gamma}C_{n}=\CC \oplus R^{\gamma}C_{n}(1-e)$. Since $[e,1-e]_{k}=0$ and  $rank_{R}(\CC) + rank_{R}(R^{\gamma}C_{n}(1-e))=n$, we conclude that $\CC^{\perp_{k}}=R^{\gamma}C_{n}(1-e)$. Thus, $\CC$ is a $\lambda$-constacyclic $k$-Galois LCD code. 
\end{proof}

\begin{prop}
Let $R$ be a finite commutative chain ring with identity and characteristic p, an rational prime, $C_{n}=\left\langle g, \, g^{n}=1\right\rangle$ be a cyclic group of order n and $R^{\gamma_{\lambda}}C_{n}$ the twisted group algebra of $C_{n}$ over $R$ where

\begin{center}
		$\gamma_{\lambda}(g^{j},g^{k})= \left\{
		\begin{array}{lll}
			\lambda,    & {\rm{if}} \,\, j+k\geq n\\
			1, &  {\rm{if}} \,\, j+k < n.
		\end{array}
		\right.$
\end{center}		

Suppose that $\lambda^{2}=1$ and  $\CC$ is a $\lambda$-constacyclic code  generated by an idempontent e. Then, $C \subset C^{\perp_{k}}$ if, and only if, $e(\sigma^{k}(e))^{\ast}=0$.
\end{prop}

\begin{proof}
Let $\CC$ be $\lambda$-constacyclic code such that $\lambda^{2}=1$ and $\CC=\left\langle e \right\rangle$ with $e^{2}=e$. If $C \subset C^{\perp_{k}}$, then $0=[e,e]_{k}=[1,e^{\ast}(\sigma^{k}(e))]$. Thus, $[\overline{g},\overline{g}e^{\ast}(\sigma^{k}(e))]=0$, for all $g \in C_{n}$. It implies $e^{\ast}(\sigma^{k}(e))=0$ and $e(\sigma^{k}(e))^{\ast}=(e^{\ast}(\sigma^{k}(e)))^{\ast}=0$.

On the other hand, if $e(\sigma^{k}(e))^{\ast}=0$, then \\ $[e, \sigma^{k}(e)]=0$, so $[\overline{g}e, \overline{h}e]_{k}=[\overline{g},\overline{h}e^{\ast}(\sigma^{k}(e))]=0$ for all $\overline{g}$, $\overline{h} \in C_{n}$. So, $[\alpha,\beta]_{k}=0$, for all $\alpha$ and $\beta \in \CC$ and the proof is completed.
\end{proof}


\begin{cor}
Let $\CC$ be $\lambda$-constacyclic code such that $\lambda^{2}=1$ and $\CC=\left\langle e \right\rangle$ with $e^{2}=e$. Then, $C \subset C^{\perp}$ if, and only if, $ee^{\ast}=0$.
\end{cor}

{\textit{Example:}}  Let $R=\F_{4}+u\F_{4}$ be a commutative local ring with $u^{2}=0$, where $\F_{4}=\{0,1,\omega, \omega + 1\}$ is the finite field with 4 elements, $C_{5}$ be the cyclic group of order 5 and $R^{\gamma}C_{5}$ be the twisted group ring where

\vr

\begin{center}
		$\gamma(g^{j},g^{k})= \left\{
		\begin{array}{lll}
			1+\omega \cdot u,    & {\rm{if}} \,\, j+k\geq 5\\
			1, &  {\rm{if}} \,\, j+k < 5.
		\end{array}
		\right.$
\end{center}		

So, in this case, we have $\overline{g}^{5}=1+\omega\cdot u$. Notice that $(1+\omega\cdot u)^{2}=1+\omega^{2}\cdot u^{2}=1$. Taking $\sigma(a+b\cdot u)=a^{2}+b^{2}(\omega + 1)\cdot u$, by \cite{YA} Proposition 1 for commutative case,  the mapping $\sigma$ is an automorphism of $R$ and $\sigma(1+\omega\cdot u)=1+\omega\cdot u$.

Now, consider the element $e=\overline{g}^{4}+ (1+\omega\cdot u)\overline{g}^{3}+ \overline{g}^{2} + (1+\omega\cdot u)\overline{g}+1$, thus

\vr

$e^{2}=\overline{g}^{8}+ (1+\omega\cdot u)^{2}\overline{g}^{6}+ \overline{g}^{4} + (1+\omega\cdot u)^{2}\overline{g}^{2}+1$

\vr

\hspace{0.3cm} $=(1+\omega\cdot u)\overline{g}^{3}+ (1+\omega\cdot u)\overline{g}+\overline{g}^{4}+\overline{g}^{2}+1 =e$ 

\vr

$e^{\ast}=(1+\omega\cdot u)\overline{g}+ \overline{g}^{2}+ (1+\omega\cdot u)\overline{g}^{3} + \overline{g}^{4} + 1=e$



\vr

$\sigma(e)=\sigma(1)\overline{g}^{4}+ \sigma(1+\omega\cdot u)\overline{g}^{3}+ \sigma(1)\overline{g}^{2} + \sigma(1+\omega\cdot u)\overline{g}+\sigma(1)=e$


\vr

\noindent since $\sigma(1+\omega\cdot u)=1+\omega \cdot u$. By Proposition \ref{GL}, the code $\CC=\left\langle e\right\rangle$ is $k$-Galois LCD code with respect the generalized $k$- Galois form.

Again using SageMath to compute the minimum distance, the code generated by $\vartheta(e)=\overline{g}^{4}+\overline{g}^{3}+ \overline{g}^{2} + \overline{g}+1$ over $\F_{4}$ is a [5,1,5] MDS code, so the code $\CC$ generated by $e$ is also MDS code over $R$. Also, the $k$-Galois dual code $\CC^{\perp_{k}}$ is a [5,4,2] MDS over $R$.  

 \section*{Data availability.}

All data are available from the authors upon reasonable request.

\section*{Conflict of interest}

All authors have participated in (a) conception and design, or analysis and interpretation of the data; (b) drafting the article or revising it critically for important intellectual content; and (c) approval of the final version.  

This manuscript has not been submitted to, nor is under review at, another journal or other publishing venue.

The authors have no affiliation with any organization with a direct or indirect financial interest in the subject matter discussed in the manuscript

\end{document}